\documentclass[12pt]{article}
\usepackage{amsthm,amssymb,amsmath,amsfonts}
\usepackage{graphicx}
\usepackage{hhline,array}
\usepackage{enumerate}
\usepackage{color}
\numberwithin{equation}{section} \numberwithin{figure}{section}

\newtheorem{theorem}{\bf Theorem}[section]

\newtheorem{Remark}{\bf Remark}
\newtheorem{example}{\bf Example}[section]

\begin{document}
\sloppy
\begin{center}
\textsc{\Large Numerical Solution of Weakly Regular Volterra Integral Equations of the First Kind}

\textbf{Denis Sidorov, Aleksandr Tynda and Ildar Muftahov}

\end{center}


\begin{abstract}

The numerical method for solution of the weakly regular scalar Volterra integral equation of the
1st kind is proposed. The kernels of such equations have jump discontinuities on the
continuous curves which starts at the origin.
The mid-rectangular quadrature rule is employed for the numerical method construction.
The accuracy of proposed numerical method is $\mathcal{O}(1/N).$

\end{abstract}

\section{Introduction}

This article deals with the following linear weakly regular Volterra integral equation (VIE) of the first kind
\begin{equation}
\int\limits_{0}^{t} K(t,s) x(s) ds = f(t), \,\,\,  0 \leq s\leq t \leq T,\,\, f(0)=0,
\label{eqi1}
\end{equation}
where kernel is defined as follows
\begin{equation*}
    K(t,s) := \left\{ \begin{array}{ll}
         \mbox{$K_1(t,s), \,\, t,s \in m_1,$} \\
         \mbox{\,\, \dots \,\, \dots \dots \dots } \\
         \mbox{$K_n(t,s), \,\, t,s \in m_n,$} \\
        \end{array} \right. \,\,  \begin{array}{ll} m_i: = \{ t, s\,\,  \bigl |\,\, \alpha_{i-1}(t) < s < \alpha_i(t) \},\\
$\,$\\
 { \alpha_0(t)=0,\,\, \alpha_n(t)=t,\, i=\overline{1,n,}}\\
\end{array}
\end{equation*}
$\alpha_i(t),$ $f(t) \in \mathcal{C}_{[0,T]}^1,$ $K_i(t,s)$
have continuous derivatives (w.r.t. $t$) for $t,s \in \overline{m_i},$  $ K_n(t,t) \neq 0,$   $\alpha_i(0)=0,$ $\,\,\,\,0 < \alpha_1(t)<\alpha_2(t)< \cdots < \alpha_{n-1}(t)<t, $
  $\alpha_1(t), \dots , \alpha_{n-1}(t) $  increase at least in the small neighborhood
$0\leq t \leq \tau,$
$0< \alpha_1^{\prime}(0) \leq \cdots \leq \alpha_{n-1}^{\prime}(0)<1. $


Such integral equations are in the core of many mathematical models in
physics, economics and ecology. The theory of integral models of evolving systems was initiated  in the early works of  L. Kantorovich, V. Glushkov
and R. Solow in the mid-20th Century.  Such theory  employs  the VIEs of the first kind where bounds of the integration interval can be  functions of time.
It is to be noted that conventional Glushkov integral model of evolving systems is the special case of the VIE (\ref{eqi1})  where all the functions $K_i(t,s)$
 are zeros except of $K_n(t,s)$.


We stress here that the VIE \eqref{eqi1} is the ill-posed problem. Such weakly regular
equations have been introduced in  \cite{Sidorov-FR}.
It is to be noted that solution of the equation
(\ref{eqi1}) may contain some arbitrary constants and can be unbounded as  $t$ goes to $0.$
Indeed, if
\begin{equation}
    K(t,s) = \left\{ \begin{array}{ll}
         \mbox{$1, \, 0< s < t/2 $}, \\
         \mbox{$-1, \, t/2 < s < t $}, \\
        \end{array} \right.
        \end{equation}
        $f(t)=t,$
        then equation (\ref{eqi1})  has the solution
         $x(t) = c -\frac{\ln t}{\ln 2},$ where $c$ reamains free parameter.
 Numerical solution of the VIE \eqref{eqi1} based on combinations of the left and right rectangle rules has been discussed by E.V. Markova and D.N. Sidorov in \cite{SidorovMarkovaAiT}.

        In this paper for such VIE with jump discontinuous kernels we propose numerical method and  discuss the analytical algorithm for construction of the continuous
solutions  in the following form:
        \begin{equation}
        x(t) = \sum\limits_{i=0}^N x_i(\ln t) t^i + t^N u(t).
        \label{eqi3}
        \end{equation}
Coefficients $x_i(\ln t)$
are constructed as polynomials on powers of $\ln t$ and they may depend on certain number of
arbitrary constants.  $N$ defines the necessary smoothness of the functions $K_i(t,s), \, f(t).$

Let us make the following notation
$
D(t) := \sum\limits_{i=1}^{n-1} \bigl| \;\alpha_i^{\prime}(t)K_n^{-1}(t,t) \bigr|\cdot\bigl |K_i(t,\alpha_i(t)) - K_{i+1}(t,\alpha_i(t))\bigr |
$
and  brifly outline the main results. Here readers may refer to the papers \cite{SidorovDU,SidorovRM}.

\begin{theorem} \label{th1}{(Sufficient Conditions of  Existence \& Uniqueness of  Local Solution)}
 Let for $t\in[0,T]$ the following conditions be satisfied: continuous $K_i(t,s), \, i=\overline{1,n}$,  $\alpha_{i}(t)$ and $f(t)$
have continuous derivatives wrt $t$,  $K_n(t,t)\neq 0$, $0=\alpha_0(t)<\alpha_1(t)< \dots < \alpha_{n-1}(t)<\alpha_n(t)=t$
 for $t \in (0,T]$, $\alpha_i(0)=0$, $f(0)=0, {D(0)<1},$ then $\exists \tau>0$ such as
eq. \eqref{eqi1} has a unique local solution in $\mathcal{C}_{[0,\tau]}.$
 Moreover if
{$\min\limits_{\tau \leq t \leq T} (t- \alpha_{n-1}(t)) = h>0$}. Then eq. \eqref{eqi1}
has {unique global solution} in  ${\mathcal{C}_{[0,T]}}$.
\end{theorem}

Let us outline the following conditions.\\
{\bf A.} Exists polynomial $\mathcal{K}_i^M(t,s) = \sum\limits_{\nu+\mu=0}^M K_{i\nu \mu} t^{\nu} s^{\mu}, i=\overline{1,n},\,\,$,
$f^M(t) = \sum\limits_{\nu =1}^M f_{\nu} t^{\nu}$, $\alpha_i^M(t) = \sum\limits_{\nu =1}^M \alpha_{i \nu} t ^{\nu}, \, i=\overline{1,n-1}$,
where $0<\alpha_{11}<\alpha_{12}< \dots < \alpha_{n-1,n}<1$
such as for $t \rightarrow +0,$ $s\rightarrow +0$
the following estimates  hold:  $|K_i(t,s)-\mathcal{K}_i^M(t,s)| = \mathcal{O}((t+s)^{M+1})$, $i=\overline{1,n}$,
$|f(t)-f^M(t)| = \mathcal{O}(t^{M+1})$,
$|\alpha_i(t) - \alpha_i^M(t)| = \mathcal{O}(t^{M+1}), i=\overline{1,n-1}$.

\noindent
{\bf B.}\, For fixed $q \in (0,1)$, $\exists \tau \in (0,T], 0<\varepsilon <1:$
$
{\max\limits_{t \in [0,\tau]} \varepsilon^{M} D(t) \leq q <1}.
$

\noindent
{\bf C}. Exists $N^*$ such as $\lim\limits_{t\rightarrow 0}  \frac{\left(\int\limits_0^t K(t,s) {\widehat{x}(s)}\, ds -f(t)\right)^{\prime}}
{t^{N^*}}  =0.$

\begin{theorem}\label{th1_2}{(Regularization)}
Let  ${\widehat{x}(t)}$ be the known function such as the condition
 {\bf C} is true for $N^* \geq M$.
  Then eq. (1)  has the  solution
${x(t) = {\widehat{x}(t)} + t^{N^*}{u(t)}},
 $ where ${u(t)} \in \mathcal{C}_{[0,T]}$ is unique and
 can be  constructed by means of the successive approximations method from the equation
$\int_0^t K(t,s) s^N u(s) \, ds = g(t),$ where $g(t):=-\int_0^t K(t,s)\hat{x(s)}\, ds +f(t).$
\end{theorem}

\begin{theorem}\label{th2}{(Relaxed Sufficient Condition for Existence \& Uniqueness)}
Let conditions {\bf B} and {\bf C} be satisfied, and
 $B(j) := K_n(0,0) + \sum\limits_{i=1}^{n-1} (\alpha_i^{\prime}(0))^{1+j} (K_i(0,0)-K_{i+1}(0,0)) \neq 0$ for $j\in {\mathbb{N}}\cup\{0 \}$.
Then eq. \eqref{eqi1} has unique  solution $x(t)=x^M(t)+t^{N^*}u(t)$ in $\mathcal{C}_{[0,T]}, \,\,$ $M \geq N.$
Moreover, for $t \rightarrow +0$ polynomial $\widehat{x}(t) \equiv x^M(t) = \sum\limits_{i=0}^M x_i t^i$ is
an $M$th order asymptotic approximation of such solution.
\end{theorem}


The paper is organized  as follows. In section 2 we propose the numerical method
for solution of the VIE \eqref{eqi1}. In section 3 we demonstrate the efficiency of proposed
numerical method on synthetic data. As footnote,  we outline the final conclusions in section 4.

\section{Numerical method}
%

In this section we propose the generic numerical method for weakly singular Volterra
integral equations \eqref{eqi1} based on the mid-rectangular quadrature rule.
The accuracy of proposed numerical method is $\mathcal{O}(1/N).$
Section 3 illustrates concepts and results of proposed numerical method on synthetic data.

For numerical solution of the equation \eqref{eqi1} on the interval \([0,T]\)
we introduce the following mesh  (the mesh can be non-uniform)
\begin{equation}\label{e2014-32}
  0=t_0<t_1<t_2<\ldots<t_N=T, \;\: h=\max\limits_{i=\overline{1,N}}(t_i-t_{i-1})=\mathcal{O}(N^{-1}).
\end{equation}
Let us search for the approximate solution of the equation \eqref{eqi1} as follows
\begin{equation}\label{e2014-33}
   x_N(t)=\sum\limits_{i=1}^{N}x_i\delta_i(t),\; t\in(0,T], \;
   \delta_i(t)=\left\{
            \begin{array}{ll}
              1, & \hbox{for } t\in \Delta_i=(t_{i-1},t_i]; \\
              0, & \hbox{for } t\notin\Delta_i
            \end{array}
          \right.
\end{equation}
with  coefficients  \(x_i,\;i=\overline{1,N}\) are under determination.

In order to find \(x_0=x(0)\) we differentiate both sides of the equation \eqref{eqi1} wrt \(t\):
\[
 f'(t)=\sum\limits_{i=1}^n \Biggl(\,\,\,\int\limits_{\alpha_{i-1}(t)}^{\alpha_i(t)}
       \frac{\partial K_i(t,s)}{\partial t}x(s)ds+\alpha'_i(t)K_i(t,\alpha_i(t))x(\alpha_i(t))-\]
     \[ - \alpha'_{i-1}(t)K_i(t,\alpha_{i-1}(t))x(\alpha_{i-1}(t))\Biggr).
\]
Therefore
\begin{equation}\label{e2014-34}
   x_0=\frac{f'(0)}{\sum\limits_{i=1}^n K_i(0,0)\left[\alpha'_i(0)-\alpha'_{i-1}(0)\right]}.
\end{equation}
Here we assume that conditions of the Theorem  \ref{th1} are satisfied and  ${\sum\limits_{i=1}^n K_i(0,0)\left[\alpha'_i(0)-\alpha'_{i-1}(0)\right]} \neq 0.$

Let's make the notation \(f_k:=f(t_k), \;k=1,\ldots,N\).
In order to define the coefficient \(x_1\) we rewrite the equation in \(t=t_1\):
\begin{equation}\label{e2014-35}
  \sum\limits_{i=1}^n \int\limits_{\alpha_{i-1}(t_1)}^{\alpha_i(t_1)}K_i(t_1,s)x(s)ds=f_1.
\end{equation}
It is to be noted that the lengths of all the segments of integration  \(\alpha_i(t_1)-\alpha_{i-1}(t_1)\) in \eqref{e2014-35} are less or equal to \(h\) and an approximate
solution is  \(x_1\) then application of the
mid-rectangular quadrature rule yields
\begin{equation}\label{e2014-36}
   x_1=\frac{f_1}{\sum\limits_{i=1}^n (\alpha_i(t_1)-\alpha_{i-1}(t_1))K_i(t_1,\frac{\alpha_i(t_1)+\alpha_{i-1}(t_1)}{2})}.
\end{equation}

 The mesh point of the mesh \eqref{e2014-32} which coincide with
  \(\alpha_i(t_j)\)  we denote as \(v_{ij}\), i.e.  \(\alpha_{i}(t_j)\in\Delta_{v_{ij}}\). Obviously \(v_{ij}<j\) for \(i=\overline{0,n-1}\), \(j=\overline{1,N}\).
{It is to be noted that \(\alpha_i(t_j)\) are not always
coincide  with any mesh point. Here \(v_{ij}\)
is used as index of the segment
\(\Delta_{v_{ij}}\),  such as  \(\alpha_i(t_j) \in \Delta_{v_{ij}}\) (or its right-hand side).}

Let us now assume the coefficients \(x_0,x_1,\ldots,x_{k-1}\) be known.
Equation \eqref{eqi1} defined in \(t=t_k\) as
\[
  \sum\limits_{i=1}^n \int\limits_{\alpha_{i-1}(t_k)}^{\alpha_i(t_k)}K_i(t_k,s)x(s)ds=f_k,
\]
we can rewrite as follows:
$I_1(t_k)+I_2(t_k)+\cdots + I_n(t_k) = f_k, $
where
$$I_1(t_k):= \sum\limits_{j=1}^{v_{1,k}-1} \int\limits_{t_{j-1}}^{t_j} K_1(t_k,s) x(s)\,ds+
\int\limits_{t_{v_{1,k}}-1}^{\alpha_1(t_k)} K_1(t_k,s) x(s)\,ds,$$
$$I_n(t_k):= \int\limits_{\alpha_{n-1}(t_k)}^{t_{v_{n-1,k}}} K_n(t_k,s) x(s)\,ds +
\sum\limits_{j=v_{n-1,k}+1}^k \int\limits_{t_{j-1}}^{t_j} K_n(t_k,s) x(s)\,ds.$$
\begin{enumerate}
\item If $v_{p-1,k} \neq v_{p,k},\, p=2,\dots, n-1$ then
$$I_p(t_k) := \int\limits_{\alpha_{p-1}(t_k)}^{t_{v_{p-1,k}}}K_p(t_k,s)x(s)\, ds +
\sum\limits_{j=v_{p-1,k}+1}^{v_{p,k}-1} \int\limits_{t_{j-1}}^{t_j} K_p(t_k,s) x(s)\,ds+$$
$$+\int\limits_{t_{v_{p,k}}-1}^{\alpha_p(t_k)} K_p(t_k,s) x(s)\, ds. $$
\item If $v_{p-1,k}=v_{p,k}$  then
$$I_p(t_k):= \int\limits_{\alpha_{p-1}(t_k)}^{\alpha_p(t_k)} K_p(t_k,s) x(s)\,ds.$$
\end{enumerate}

%
%
\begin{Remark}\label{rem1}
The number of terms in each line of the last formula depends on an array \(v_{ij}\),
defined using the input data: functions \(\alpha_i(t),\;i=\overline{1,n-1}\),
and fixed (for specific $N$) mesh. 
\end{Remark}
Each integral term we approximate using the mid-rectangular quadrature rule, e.g.
\[
  \int\limits_{t_{v_{p,k}-1}}^{\alpha_p(t_k)}K_p(t_k,s)x(s)ds\approx
  \left(\alpha_p(t_k)-t_{v_{p,k}-1}\right)K_p\left(t_k,\frac{\alpha_p(t_k)+t_{v_{p,k}-1}}{2}\right)
  x_N\left(\frac{\alpha_p(t_k)+t_{v_{p,k}-1}}{2}\right).
\]
Moreover,
on those intervals where the desired function has been already determined,
we select \(x_N(t)\) (i.e. \(t\leqslant t_{k-1}\)).

On the rest of the intervales an unknown value $x_k$ appears in the last terms.
We explicitly define it and proceed in the loop for \(k \). The number of these terms is determined from the initial data \(v_{ij}\) analysis.
The accuracy of the numerical method is  \(O\left(\frac{1}{N}\right)\).



\section{Numerical illustrations}

Let us consider three examples. In all cases the uniform mesh is used.

\begin{example} \label{ms_ex1} 
\begin{equation*}
\int\limits_0^{t/3} (1+t -s) x(s)\,ds-\int\limits_{t/3}^{t} x(s)\,ds =\frac{t^4}{108} - \frac{25 t^3}{81},\; t\in [0,\,2],
\end{equation*}
exact solution is $\bar{x}(t)=t^2$.
Tab. \ref{ms_tab1} demonstrates the errors   $\varepsilon=\max\limits_{0\leq i \leq n}|\bar{x}(t_i)-x^h(t_i)|$ for various steps $h$.

\begin{table}[htb]\label{ms_tab1}
\caption{Errors for the 1st Example.}
\begin{center}
\begin{tabular}{c|l}
\hline
 ${h}$ & $\varepsilon $   \\  \hline
 $1/32$   &  0{.}13034091293670258 \\
 $1/64$   &  0{.}07804538180930365 \\
 $1/128$  &  0{.}03989003750757547 \\
 $1/256$  &  0{.}01975354947865071 \\
$1/512$  &  0{.}010027923872257816\\
$1/1024$  &  0{.}005083865773485741\\
$1/2048$  &  0{.}0025693182974464435\\
$1/4096$  &  0{.}001288983987251413\\
$1/8192$  &  0{.}0006500302042695694\\ 
\end{tabular}
\end{center}
\end{table}

\end{example}

\begin{example}\label{ms_ex2}

\begin{equation*}
\int\limits_0^{\frac{t}{9}}(1+t-s)x(s)\,ds-\int\limits_{\frac{t}{9}}^{\frac{2t}{9}}x(s)\,ds
-2\int\limits_{\frac{2t}{9}}^{\frac{4t}{9}} x(s)\, ds
 +
\int\limits_{\frac{4t}{9}}^t x(s)\, ds=\frac{11t^4}{26244}+\frac{547t^3}{2187},\; t\in [0,\, 2],
\end{equation*}
exact solution is $\bar{x}(t)=t^2$.
Tab. \ref{ms_tab2} demonstrates the errors $\varepsilon$  for various steps $h$.

\begin{table}[htb]\label{ms_tab2}
\caption{Errors for the 2nd Example.}
\begin{center}
\begin{tabular}{c|l}
\hline
 ${h}$ & $\varepsilon$   \\  \hline
 $1/32$   &  0{.}13718808476353672 \\
 $1/64$   &  0{.}07408554651043886 \\
 $1/128$  &  0{.}04531351578371812 \\
 $1/256$  &  0{.}022111520501482573 \\
$1/512$  &  0{.}011079518173630731\\
$1/1024$  &  0{.}005492567505257284\\
$1/2048$  &  0{.}0027453216364392574\\
$1/4096$  &  0{.}0014125244842944085\\
$1/8192$  &  0{.}00077170109943836\\ 
\end{tabular}
\end{center}
\end{table}
\end{example}

\begin{example}\label{ms_ex3}

\begin{equation*}
2\int\limits_0^{\sin\frac{t}{2}} x(s)\,ds-\int\limits_{\sin\frac{t}{2}}^{2\sin\frac{t}{3}}x(s)\,ds +
\int\limits_{2\sin\frac{t}{3}}^t x(s)\, ds=\frac{t^3}{3}+\sin^3\frac{t}{2} - \frac{16}{3}\sin^3\frac{t}{3},\; t\in \left[0,\, \frac{3\pi}{2}\right],
\end{equation*}
exact solution is $\bar{x}(t)=t^2$.
Tab. \ref{ms_tab3} demonstrates the errors $\varepsilon$  for various steps $h$.

\begin{table}[htb]\label{ms_tab3}
\caption{Errors for the 3rd Example.}
\begin{center}
\begin{tabular}{c|l}
\hline
 ${h}$ & $\varepsilon$   \\  \hline
 $1/32$   &  1{.}2810138805937967 \\
 $1/64$   &  0{.}7064105257311724 \\
 $1/128$  &  0{.}3172969521937503 \\
 $1/256$  &  0{.}16990268475221626 \\
$1/512$  &  0{.}11787087222029413\\
$1/1024$  &  0{.}07940422358498633\\
$1/2048$  &  0{.}06518995509284764\\
$1/4096$  &  0{.}06004828109245386\\
$1/8192$  &  0{.}046102790104048275\\ 
\end{tabular}
\end{center}
\end{table}
\end{example}

\section{Conclusion}

In this article we addressed the novel class of weakly regular linear Volterra integral
equations of the first kind first introduced in \cite{Sidorov-FR}. We outlined the main
results for this class of equation previously derived. The main contribution
of this paper is a generic numerical method designed for solution of such weakly regular equation.
The numerical method employes the mid-point quadrature rule and enjoy the
the $\mathcal{O}(1/N)$ order of accuracy.  The illustrative examples demonstrate
the efficiency of proposed method.  As footnote let us outline that proposed
approach enable construction of the 2nd order accurate numerical method.
This improvement will be done in our further works.



\begin{thebibliography}{99}

\bibitem{SidorovDU}
D.N. Sidorov,
On parametric families of solutions of Volterra integral equations of the first kind with piecewise smooth kernel, \emph{Differential Equations,} Vol. 49, N.2, 2013, 210--216.

\bibitem{SidorovRM}
{D.N. Sidorov}, { Solution to {s}ystems of {V}olterra {i}ntegral
  {e}quations of the {f}irst {k}ind with {p}iecewise {c}ontinuous {k}ernels},
  \emph{Russian Mathematics}, Vol. 57, 2013, 62--72.

\bibitem{SidorovMarkovaAiT}
E.V. Markova, D.N. Sidorov,
On one integral Volterra model of developing dynamical systems,
\emph{Automation and Remote Control,}
 Vol. 75, No. 3, 2014, 413--421.



  \bibitem{Boikov-Tynda11}
  Boikov I.V., Tynda A.N., Approximate solution of nonlinear integral equations of
  developing systems theory. {\em Differential Equations}, Vol.39, 9, 2003, 1214-1223.

  \bibitem{Tynda-18}
   A.N. Tynda,  Numerical  algorithms of optimal complexity for weakly singular
     Volterra integral equations, {\em Comp. Meth. Appl. Math.},
     Vol.6 (2006) No. 4, 436--442.
  \bibitem{Tynda-pamm1}
   A.N. Tynda,  Numerical methods for 2D weakly singular Volterra integral equations
                of the second kind. \emph{PAMM}, Volume 7 (2007), Issue 1.
\bibitem{Sidorov-FR}
D.N. Sidorov, Volterra Equations of the First kind with
Discontinuous Kernels in the
Theory of Evolving Systems Control. \emph{Stud. Inform. Univ.}, Vol. 9 (2011), 135--146







\end{thebibliography}
\end{document}